**Towards Industry 5.0: A Systematic Literature Review on Sustainable and Green Composite Materials Supply Chains.**

Md Rabiul Hasan[1[0000-0002-5200-253X]], Muztoba Ahmed Khan[2[0000-0002-2668-1057]], Thorsten Wuest[1*[0000-0001-7457-7927]]

[1]IMSE, West Virginia University, Morgantown, WV 26501
[2]The Sullivan School of Business and Technology, Carroll University, Waukesha, WI 53186
mh00071@mix.wvu.edu; mkhan@carrollu.edu;thwuest@mail.wvu.edu

* corresponding author: thwuest@mail.wvu.edu

## Abstract

Sustainable supply chain management is a key objective of Industry 5.0, leveraging technologies like real-time data analytics, connectivity, and intelligent automation. At the same time, composite materials present benefits such as lightweight structures, crucial for reducing fuel consumption. This study investigates the intersection between sustainable supply chains and composites by analyzing the current status, research gaps, methodologies, and future research opportunities through bibliometric analysis and a systematic review of the state of the art in the composite materials supply chain. A systematic literature review approach is employed to analyze the Scopus and Web of Science (WOS) databases, offering a comprehensive overview of the existing literature. Through bibliometric analysis, the study investigates countries, authors, citations, keywords, subject areas, and article themes within the metadata to provide additional context. An in-depth analysis of thirty selected papers (n=30) sheds light on the key contributions, major challenges, and Key Performance Indicators (KPIs) across various instances of composite material supply chains, resulting in a generalized overview. Furthermore, this research suggests future directions to link the sustainability efforts in composite materials supply chains with current research gaps. The study underscores diverse research themes in the field, highlighting a few influential works and presenting opportunities for advancement in this emerging area. Collectively, these findings offer valuable insights and a robust roadmap for future research in this domain.

## 1. Introduction

In today's world, environmental consciousness for the protection of the environment has been seen among different organizations. To establish sustainable industries, contemporary organizations are looking for a green journey (Bansal and Roth, 2000). The term "Green" refers to incorporating different environmental aspects. According to a definition by Srivastava (2007), Green Supply Chain (GrSC) integrates environmental thinking starting in the supply chain product design phase and also includes the end-of-life management of the product after its useful life, like recycling of that product. The definition of GrSC varies in the literature, and is sometimes defined as 'supply chain environmental management', 'green purchasing and procurement', 'sustainable supply network management' (Tseng et al., 2019). Min and Galle (2001) defined green purchasing as a process of environment-conscious purchasing that encourages recycling and reducing waste sources. Most broadly, GrSC is sometimes also considered Sustainable Supply Chain Management (SSCM). SSCM includes, for Pagell and Shevchenko (2014), not only economic viability but also not harming social and environmental systems. Many criteria are used for environmental measures on the supply chain. It relates mainly to two groups: first, the product or service being purchased, and second the manufacturing or supplier providing them (Igarashi et al., 2013). Recently, waste management and reverse logistics have gained more attention for recycling different contemporary materials like plastics, electronics materials, organic items, glasses, etc. In the study of the waste reverse supply chain, Van Engeland et al. (2020) described it as the process of recapturing or creating value, including the typical transportation, collection, recovery, and disposal of waste. The closed-loop supply chain network optimization problem is another field for recycling management in the supply chain (Zohal and Soleimani, 2016). Supply chain optimization network design also depends on the material types recovered in a recycling process. For example, different steel products are made from steel scraps. Thus, it positively impacts the environment by reducing environmental pollution (Cui et al., 2017).

Many researchers have utilized mathematical models to reduce transportation costs and decrease environmental pollution. Wang et al. (2011) proposed a multi-objective optimization function for network design considering carbon emissions. It was suggested that improving network capacity decreases carbon emissions and total costs. In 2017, a similar study was also conducted by Nurjanni et al. (2017), who included forward and reverse networks as a closed-loop supply chain. Zhao et al. (2018) analyzed the impact of reverse logistics for electronic products by system modeling with Vensim software, a computer simulation-based cause and effect feedback process. This study confirmed that the combination of a higher recovery ratio and transparency reinforcement in GrSC has a significant impact on reducing the bullwhip effect. Again, many studies have utilized mixed-integer programming models for the closed-loop supply chain network (Özceylan and Paksoy 2013; Jindal and Sangwan 2014; Jabarzadeh et al. 2020). Several studies also focused on stochastic programming models where they tried to develop optimization in a closed-loop supply chain (Giri and Masanta 2020; Kalantari Khalil Abad and Pasandideh 2022; Yılmaz et al. 2021). In a two-stage stochastic programming model, the first stage is the strategic supply chain network design (SCND) decision. The second is the tactical SCND, where production quality and material flow in the supply chain are determined. Interestingly, in 2015, Rezaee et al. (2017) found that carbon emission pricing greatly impacts supply chain configuration. Their discrete GRSC network design with a two-stage stochastic programming model pointed out that the supply chain configuration would be greener if the carbon price was higher.

In GrSC, an integrated model is needed to provide an optimal solution embracing environmental, economic, and societal sustainability (Meng et al., 2020). A review paper by Oliveira and Machado (2021) draws attention to the lack of research on the whole spectrum of the closed-loop supply chain (economic, environmental, and social). The performance assessment of GRSC based on the suppliers' viewpoints is likely to increase. There is also an excellent opportunity to work with big data for GrSC (Tseng et al., 2019). Again, developing a model with social welfare under uncertain supply and demand in the reverse supply chain in terms of environmental benefits still needs to be solved (Ye et al., 2016).

This study aims to analyze the state of the art of the composite materials supply chain to identify research gaps, future research directions, and challenges, and Key Performance Indicators (KPIs) to support and enhance green and sustainable practices in the composite materials supply chain. This paper also explores the meta-data for bibliometric analysis to assess the impact, trends, and patterns of published scholarly articles. The rest of the paper is organized as follows: Section 2 describes the methodology and framework of this study. Section 3 represents the bibliometric analysis. We explored the comprehensive review of the selected articles in Section 4, including the significant outcomes of each paper, composite materials supply chain and recycling, development of KPI's model, challenges and future research direction. In Section 5, we conclude the study with a summary and limitations of this research work.

## 2. Method

For this study, a Systematic Literature Review (SLR) has been conducted. According to a definition provided by Okoli and Schabram (2010), the SLR is "*a systematic, explicit, and reproducible method for identifying, evaluating, and synthesizing the existing body of completed and recorded work produced by researchers, scholars, and practitioners*". SLR has been chosen for this research work due to its accurate and reliable outcomes while reviewing and synthesizing any research topics. This study's SLR follows the Preferred Reporting Items for Systematic Reviews and Meta-Analysis (PRISMA) framework. In 1999, QUOROM (Quality of Reporting of Meta-analysis) was developed for reporting guidelines. However, due to the inability to capture the complete scenario, a new framework named PRISMA was developed as an evolution of the original QUOROM guideline, which contains a 27-item checklist and a four-phase flow diagram. Although the PRISMA framework is not a quality assessment tool, it confirms the transparency and clarity of the reporting of systematic studies (Liberati et al., 2009).

### 2.1. Data Identification

To review the comprehensive published research works, both major scientific databases (WoS and Scopus) have been considered for this study. Covering major publications from major publishers, more than 13 million (Scopus) and 10 million (WoS) documents are indexed, providing a comprehensive selection (Franceschini et al., 2016).

These two databases' advanced query search with indexed papers ensures the quality of scientific research, which is missing in the Google Scholar database. The search time horizon ranges from 1977 to 25 May 2023 to cover all relevant previous publications. To identify published works in the supply chain domain, some key terms like "supply chain," "logistic," "network design," or "facility location" have been used. Furthermore, composite materials can be broadly categorized based on matrix (polymer, metal, ceramic, carbon) and fiber compositions (Daniel and Ishai, 2006) . In this study, the composites and composites material are defined as having the property of the material, which includes two or more constituent materials mostly found as a matrix or reinforced phase of the material (Harik and Wuest, 2020). Hence, different keywords have been used focusing on different composite materials and their types, such as "composite material", "recycling composite", "flake", "wood plastic", "fiberglass", "polymer matrix", "ceramic matrix", "carbon fiber". Boolean operators are used to integrate the supply chain and composite materials keywords to broaden the search area. The search terms have been cross-checked with all investigators of this research work and refined accordingly from the initial one to validate the search string and terminologies. This ensures that the search vocabulary is strictly followed while identifying the exact search radius. The following search string is used: ((("supply chain*") OR (logistic*) OR ("network design*") OR ("facility location*")) AND (("composite material*") OR ("recycling composite*") OR ("flake*") OR ("wood fiber") OR ("wood plastic") OR ("fiberglass") OR ("polymer* matrix") OR ("ceramic* matrix") OR ("Carbon fiber*") OR ("reinforced plastic*") OR ("metal* matrix") OR ("reinforced composites")))

## 2.2.Reviewing data for eligibility

In the initial search using the developed search string, the Scopus database yielded 677 results, while the WOS found 255 papers. Additionally, book chapters and business reports are excluded, and only conference and journal articles are considered. Combining both the results with the filtration of only the English language and removing duplicated articles by the Excel duplicate functions resulted in 744 unique articles. In the following step, the 744 unique articles were screened by reading the title and abstract, and any conflicts were resolved by discussing the extraction criterion. In this step, it was evident that most papers investigated the composite materials' properties instead of focusing on the composite materials' supply chain. Thus, this screening stage resulted in 101 papers for full-text review.

After the full-text review, 24 papers closely aligned with the research objectives. Following the full-text review, six additional research articles were considered by reference tracking the 24 selected eligible papers. The PRISMA flow diagram used in this study is shown in Figure 1. In the full-text screening phase, the dominant reasons for excluding research articles are:

1. Investigating the composite materials recycling process improvement rather than the supply chain
2. Focused on the composite materials quality and property development.
3. Not focusing on the supply chain perspectives and drivers
4. Not a peer-reviewed journal article, for example, a trade journal

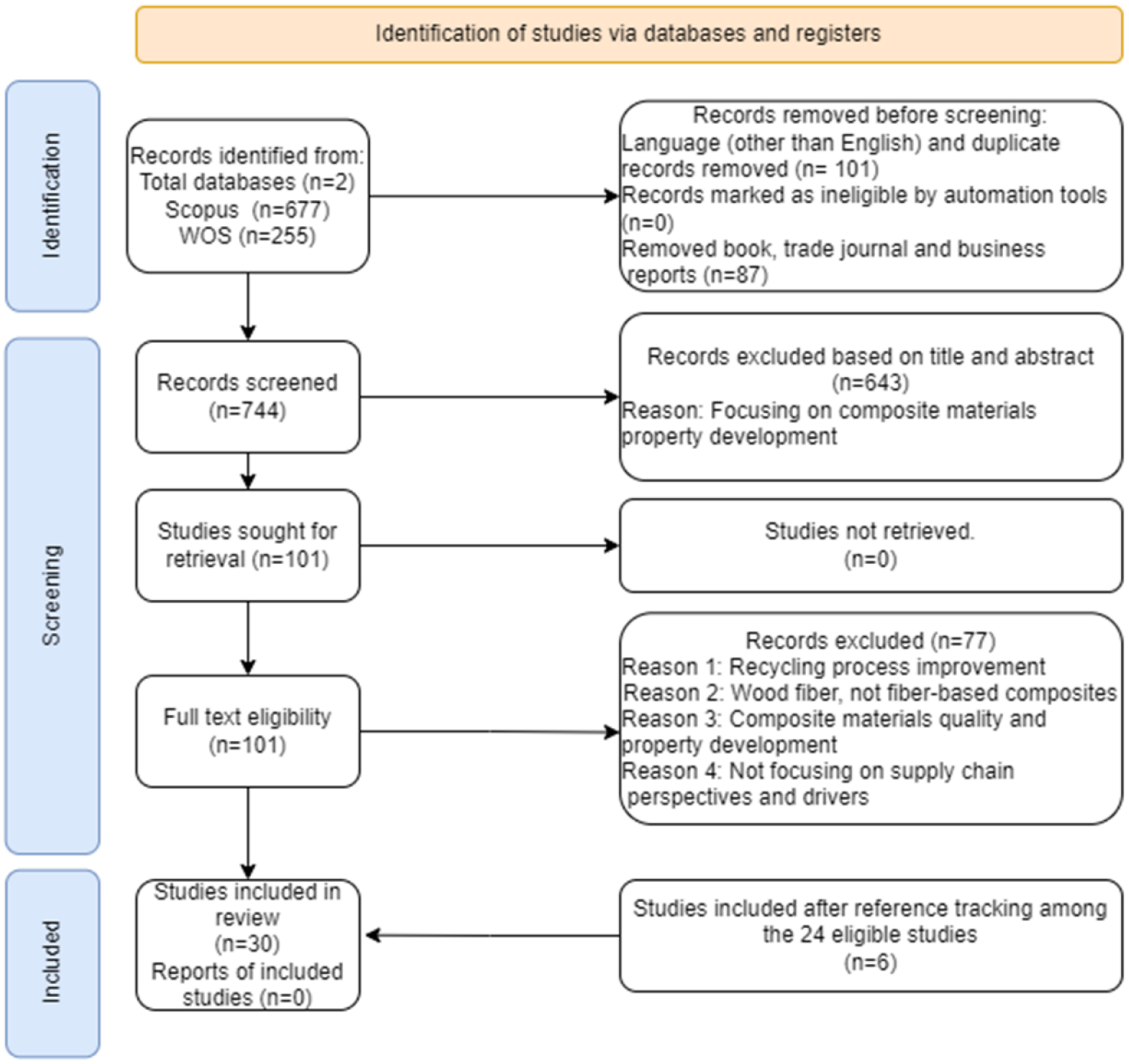

Figure 1: PRISMA Framework (Page et al., 2021)

## 3. Bibliometric analysis

For the bibliometric analysis, different research papers utilized different software packages with different limitations and capabilities. Some popular software tools are Vosviewer, Gephi, Publish or Perish, Bibexcel, Hitsite, and Bibliometrix (Biblioshiny), which is a package in R programming language. For this study, Vosviewer (v1.6.18) and Bibliometrics (Biblioshiny) libraries in R programming language have been selected due to their capability for comprehensive network analysis (Aria and Cuccurullo 2017; van Eck and Waltman 2014; Rahman et al. 2023; Khosravi et al., 2023). The bibliometric data of the chosen articles is included in CSV format combining the two databases (Scopus and WOS). The following sections describe the publication types and sources, author and citation analysis, words and related themes, and publication by countries.

### 3.1. Publication types and sources

Among the 30 papers analyzed, 46% (n=14) are conference papers, and 54% (n=16) are journal articles. In terms of publication sources, the majority of these studies appear in the Journal of Cleaner Production, followed by the Journal of Advances in Transdisciplinary Engineering, as illustrated in Figure 2. And it can be referred that publications are diverse in this research field.

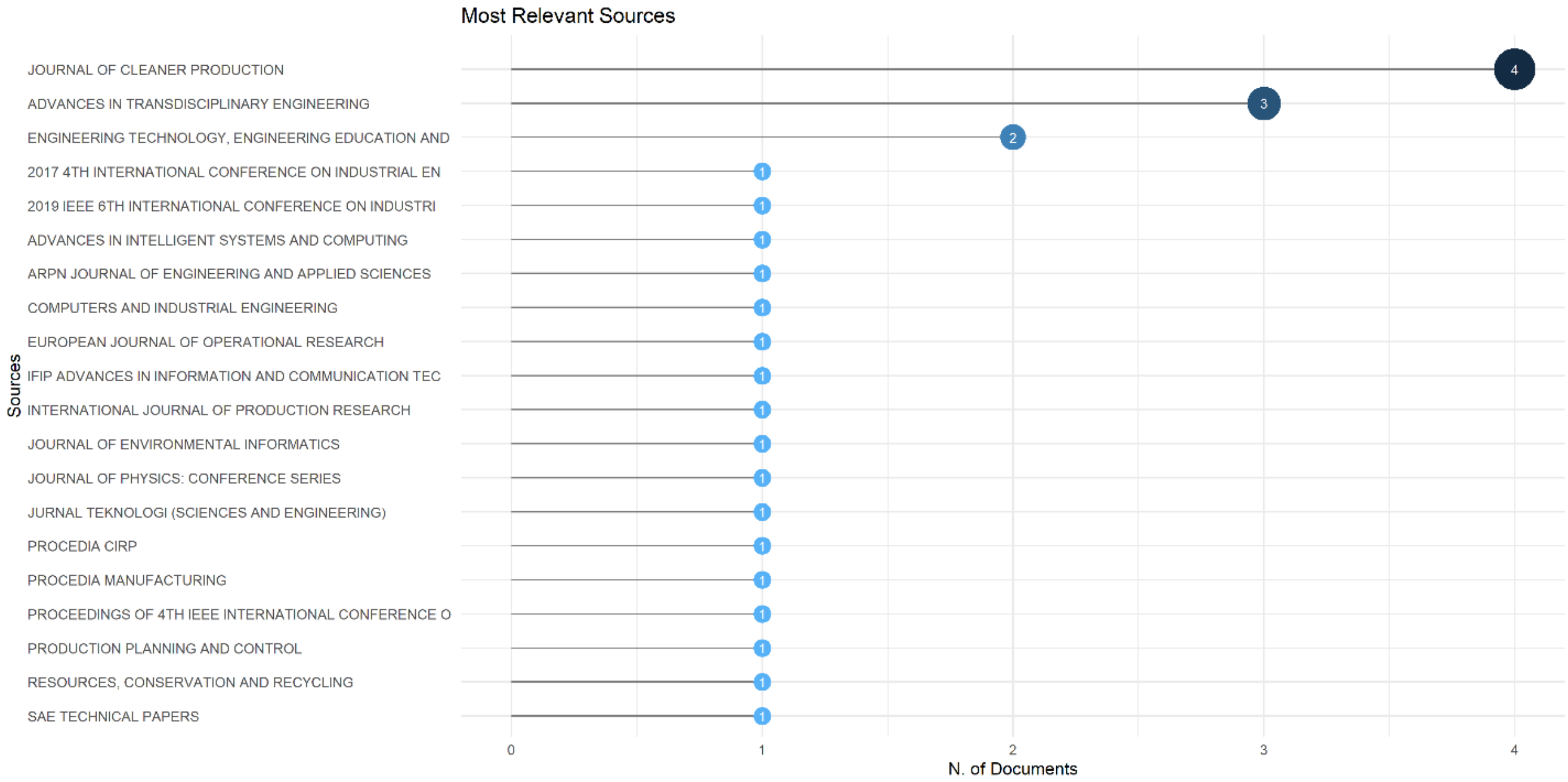

Figure 2: Sources of publication

Although these two journals dominate in terms of number of publications, the Resources, Conservation, and Recycling journal has the highest impact factor (13.2) and cite score (20.3) among the sources of publications.

### 3.2.Author and citation analysis

The author with the greatest number of papers is Mondragon A.E. Coronado (8 out of 30). In terms of number of published articles, Hogg PJ and Mastrocinque E. ranked second and third, respectively. Four other authors published three articles each, while the rest published two, based on the top ten authors by number of publications in this research field, as shown in Figure 3.

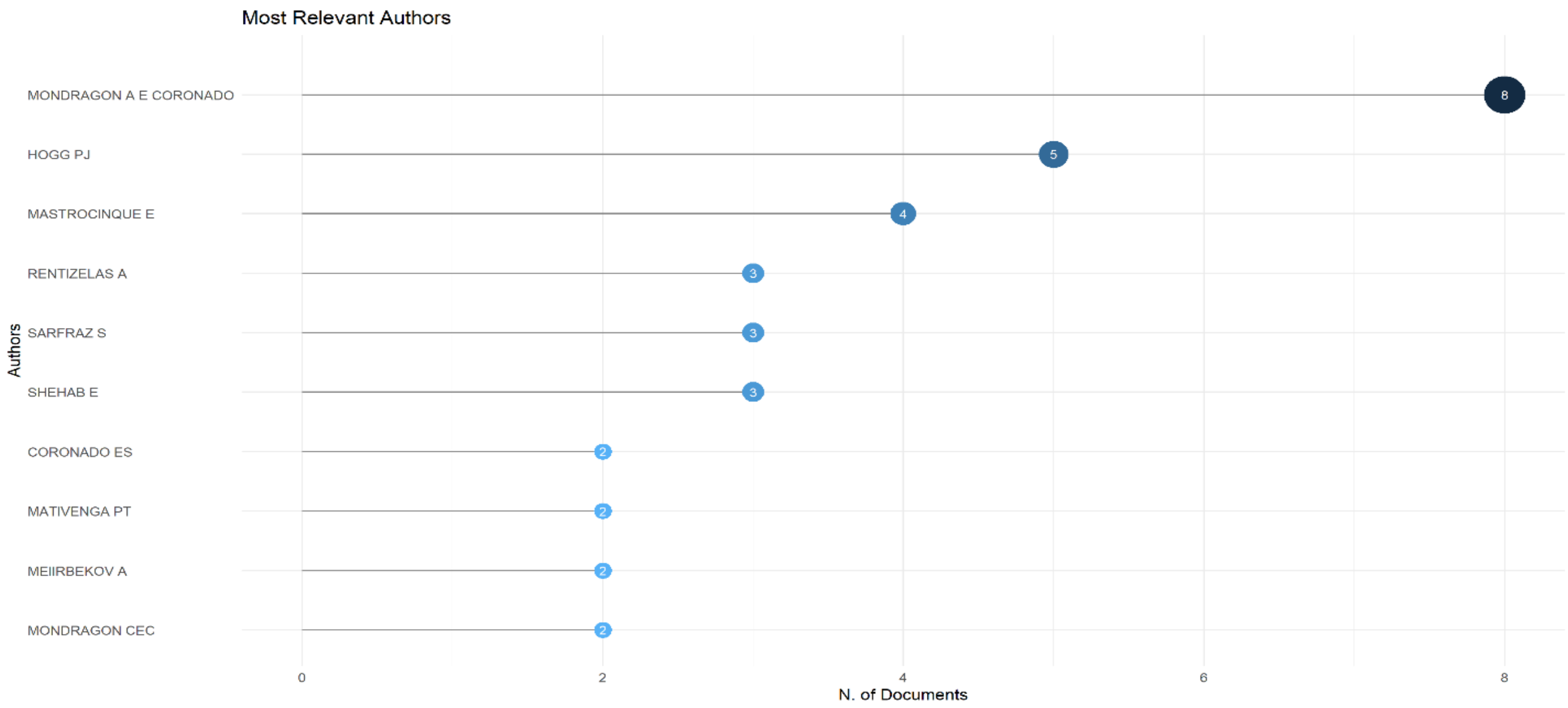

Figure 3: Number of publications of authors.

The relationship between cited and citing documents can be found in the citations (Nicolaisen, 2007). The top ten authors based on total citations (h_index, g_index, m_index) are shown in Table 1 with their year of first publication. Mondragon Adrian E. Coronado has the highest number of citations, with his first publication in 2015. Etienne S. Coronado Mondragon and Paul J. Hogg came next, ranking two and three on the list, respectively. Additionally, Adrian E. Coronado Mondragon has the highest h_index, g_index, and m_index. This bibliometric implies that Mondragon Adrian E. Coronado has had a substantial impact on the supply chain for composite materials. Furthermore, it can be inferred from the authors' first year of publication that the trends and subjects covered by this field of study are relatively recent—only a few years old. In addition, Hogg PJ has four publications, and his h-index of 2 indicates that at least two of these have received two or more citations.

Table 1: Top ten author citation

| Element | h_index | g_index | m_index | TC | NP | PY_start |
|---|---|---|---|---|---|---|
| Adrian E. Coronado Mondragon | 4 | 8 | 0.444 | 99 | 8 | 2015 |
| Etienne S. Coronado Mondragon | 2 | 2 | 0.333 | 74 | 2 | 2018 |
| Paul J. Hogg | 2 | 4 | 0.222 | 22 | 5 | 2015 |
| Ernesto Mastrocinque | 2 | 3 | 0.222 | 13 | 4 | 2015 |
| Paul Tarisai Mativenga | 2 | 2 | 0.222 | 64 | 2 | 2015 |
| Arshyn Meiirbekov | 2 | 2 | 0.500 | 6 | 2 | 2020 |
| Christian E. Coronado Mondragon | 2 | 2 | 0.286 | 63 | 2 | 2017 |
| Sarah Oswald | 2 | 2 | 1 | 26 | 2 | 2022 |
| Athanasios Rentizelas | 2 | 3 | 1 | 29 | 3 | 2022 |
| Shoaib Sarfraz | 2 | 2 | 0.4 | 8 | 3 | 2019 |

The impact of local citation serves as another name for the immediate collaboration network. The global citation, on the other hand, refers to citations that are part of the greater community and not directly from direct collaborators (Fahimnia et al., 2015). There is a gap between the authors' local and global citation impacts in the clustering by coupling based on their co-authorship pattern within the citation context as described in Table 2. This indicates that the supply chain for composite materials has drawn attention from researchers in other fields, which is not reflected in our search results. It is also interesting to find that the rank of local citation does not match the global citation. For example, it suggests that Etienne S. Coronado Mondragon is more prevalent in other research fields than composite materials supply chains, as his global citation score is far higher than his local citation score.

Table 2: Clustering of coupling- Authors (Measured by references)

| Authors | Impact measure- Local citation score | Impact measure- Global citation score |
|---|---|---|
| Ernesto Mastrocinque | 0.88 | 3.25 |
| Paul J. Hogg | 0.90 | 4.4 |
| Adrian E. Coronado Mondragon | 0.94 | 13.9 |
| Etienne S. Coronado Mondragon | 1 | 37 |
| Christian E. Coronado Mondragon | 1 | 31.5 |
| Athanasios Rentizelas | 1 | 9.7 |
| Sarah Oswald | 1 | 13 |
| Stefan Siegl | 1 | 13 |
| Nikoletta L. Trivyza | 1 | 2.5 |
| Paul Tarisai Mativenga | 1.25 | 32 |

### 3.3.Words and related themes

To gain a complete overview of a paper, the author's keywords play a significant role in representing that. Similar keywords are identified and replaced with identical keywords. For example: "Composites," "composite," or "composite materials" are replaced with the keyword "composites" only. For this study, the co-occurrence of the author keywords has been investigated. Due to the small number of articles (n=30) meeting the review criteria, the minimum number of occurrences of an author keyword is considered one, resulting in 76 unique keywords. In the network visualization of author keywords, 13 clusters have been found among the 74 keywords as shown in Figure 4.

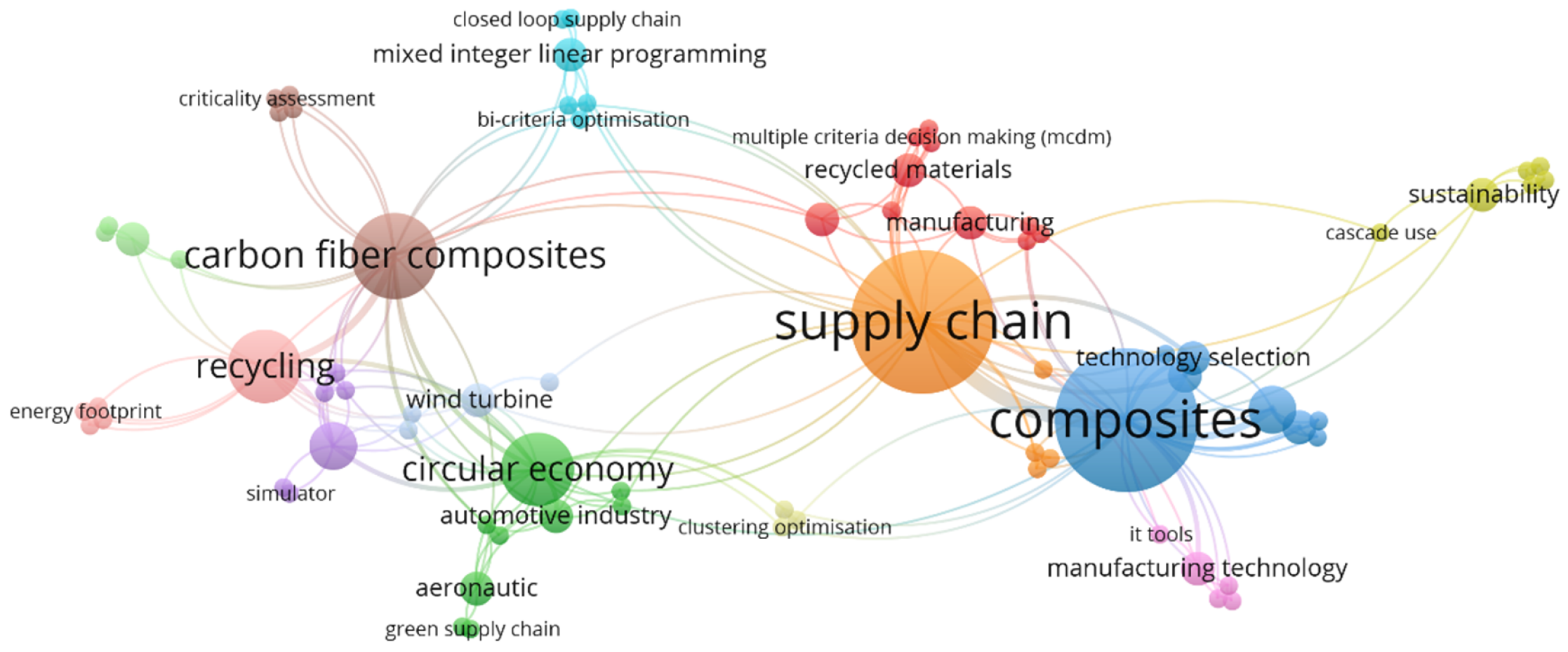

Figure 4: Co-occurrence of keywords (VosViewer)

From this analysis, it has been inferred that there is a diverse research focus in the composite manufacturing supply chain, and a link exists between different research areas. To better understand the relationship between major topics, we conducted a thematic evolution analysis of selected abstracts from articles published in four time periods: [2002-2017], [2018-2019], [2020-2021], and [2022-2023]. As shown in Figure 5 below, substantial interconnection has been found in the field of study research. The characteristics of the lines determine the significance of the connections between different words. It is apparent from the figure that the research topics are correlated with each other. More research has been carried out in the sustainability domain from 2002 through 2017. Furthermore, it can be seen that there is a trend of research ongoing in the recycling and optimization topics in the field of study.

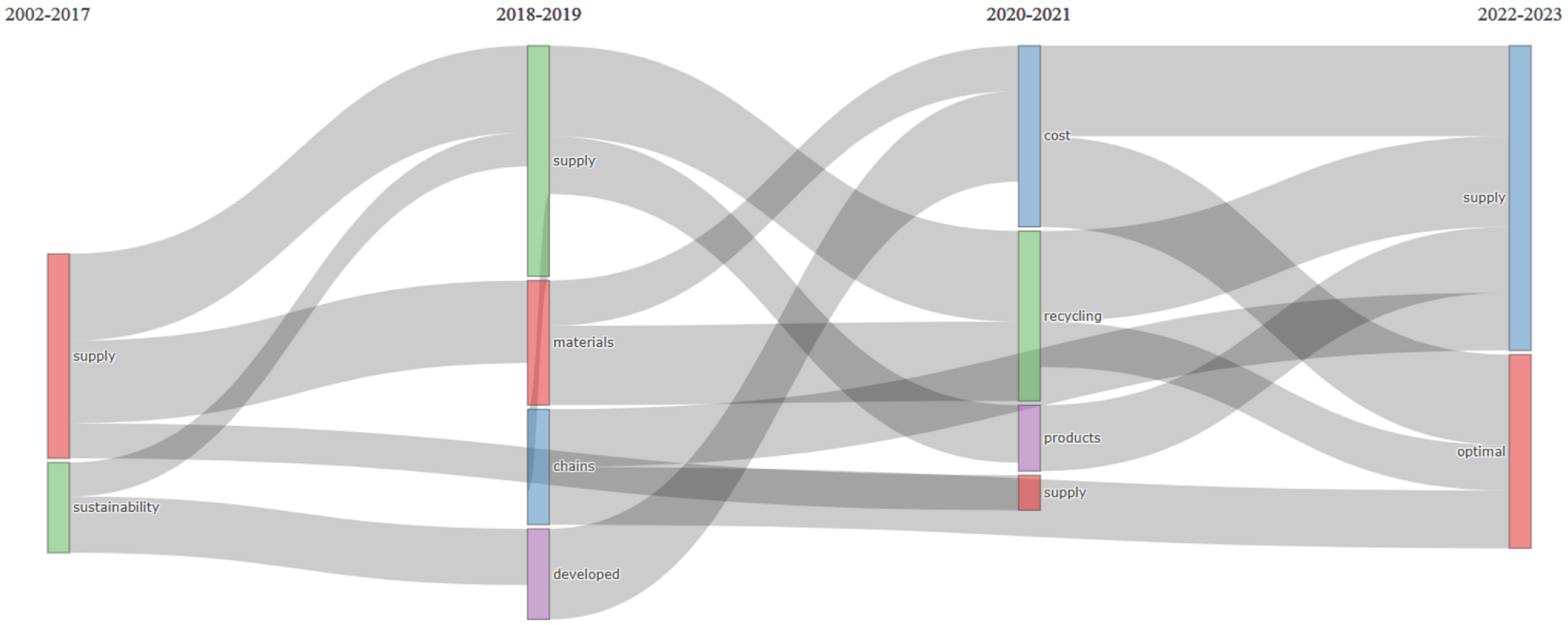

Figure 5: Thematic evaluation keywords

Another simple way to understand the prevalent themes in a complex setting is to explore the word cloud of the author's keywords. Figure 6 shows the word cloud of the referenced articles where more prominent and bolder fonts represent the frequency of that phrase, and fewer and smaller ones are vice versa.

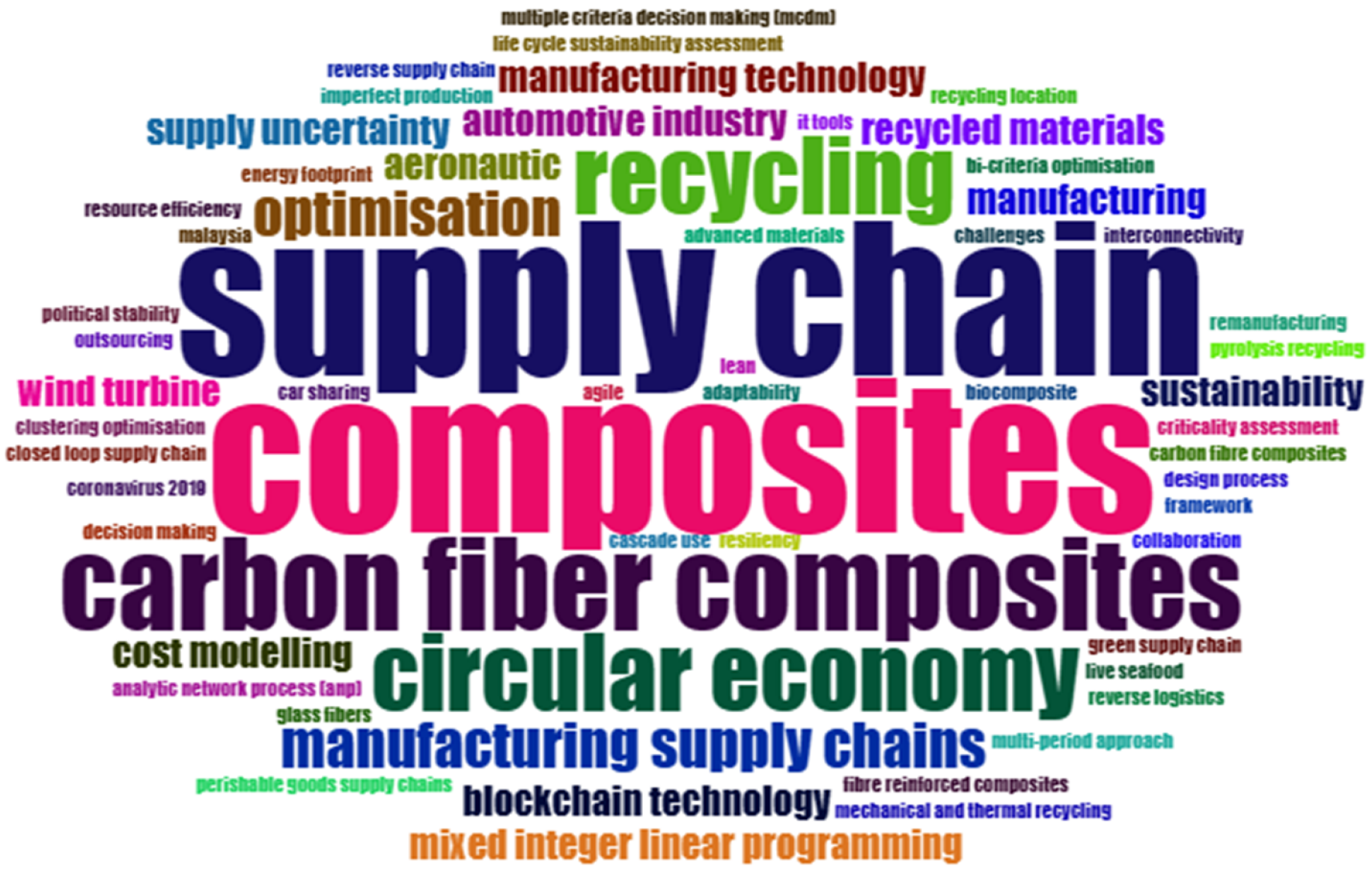

Figure 6: Keywords word cloud

After analyzing the author's keywords word cloud, some interesting themes were identified, such as blockchain technology, remanufacturing, political stability and reverse logistics. The word frequency of the keywords is shown in Figure 7. Industry-specific terms, such as automotive and aeronautics, appear frequently in the word analysis, indicating industry-specific research articles on the composite materials supply chain. Overall, the above results indicate potential interdisciplinary research for this study.

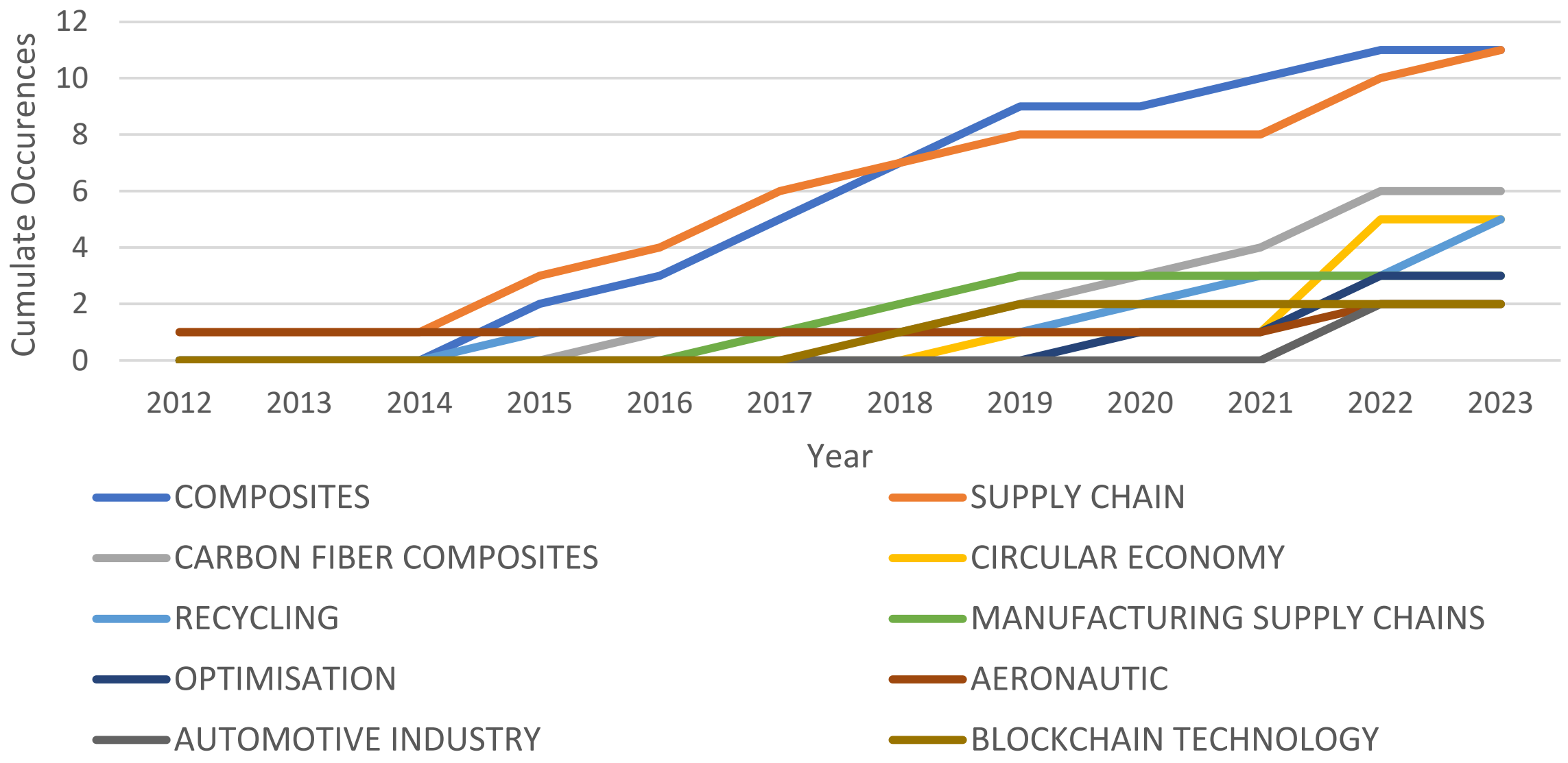

Figure 7: Word frequency of the keywords

### 3.4. Publication by countries

This section examines nations' citation counts, publication totals, and collaboration networks. Based on the data in Table 3, it has been found that the United Kingdom has the highest number of citations in this research field, with a total of 200 citations. Following the UK are Taiwan, Canada, and France, respectively, with considerable citations in the same field. It is interesting to note that these countries have made significant contributions to this area of study, and their research outputs have been widely recognized and cited by other scholars and researchers. It is surprising to find that the United States of America is not among the top ten countries based on citation number. Figure 8 illustrates the scientific contributions of various countries, with deeper blue representing more publications in the field of composite materials supply chain.

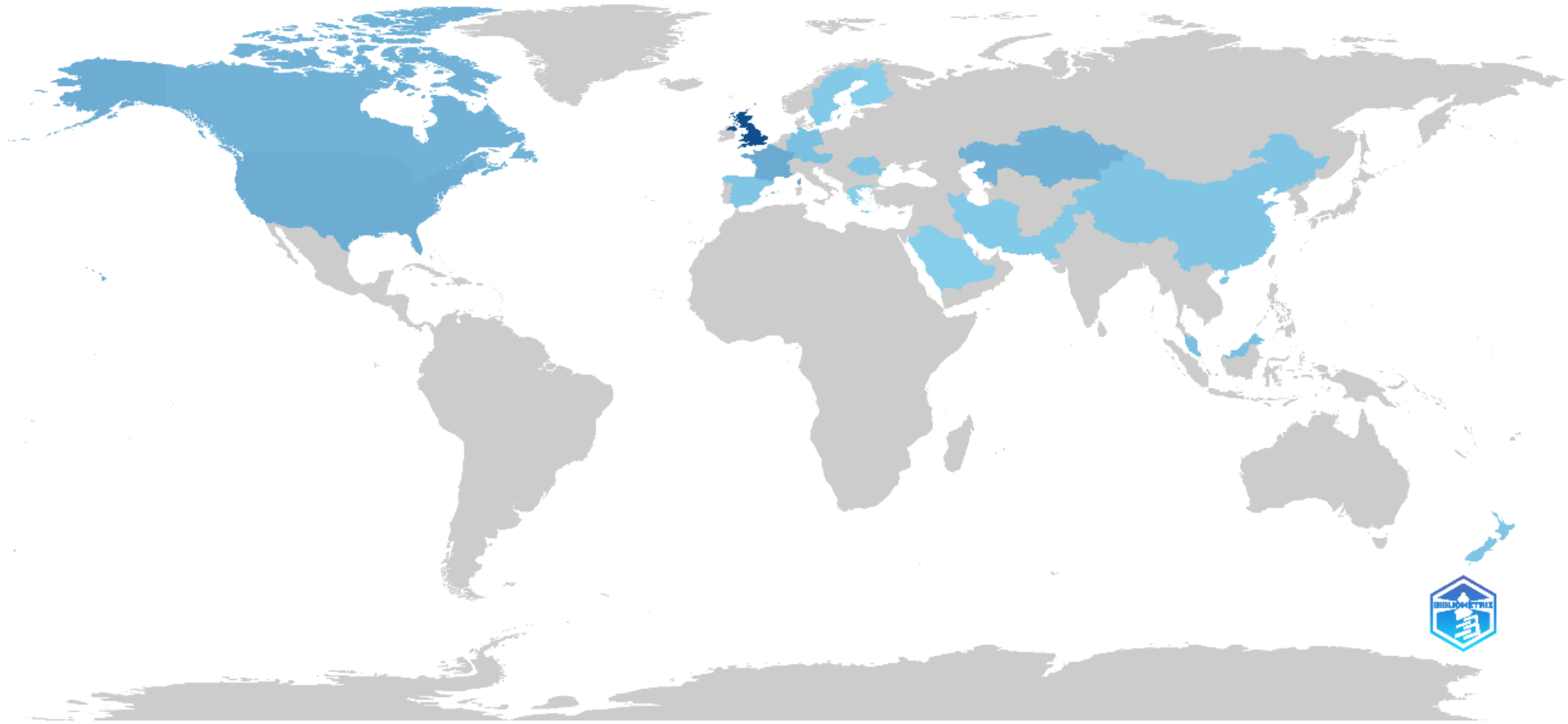

Figure 8: Map of countries' scientific output

The co-authorship of countries demonstrates collaboration between researchers from different nations. The minimum number of citations is taken into consideration by the VosViewer software for both documents and citations. Only 18 countries meet this threshold out of 20 total. The UK has the most published articles (n=15) in this field. Although the US and Greece have the same number of articles (n=4), the US does not have a significant number of citations, as shown in Table 3.

Table 3: Publications between countries in terms of the amount and citations.

| ID | Country | Documents | Citations |
|---|---|---|---|
| 1 | United Kingdom | 16 | 200 |
| 2 | Greece | 4 | 29 |
| 3 | Germany | 2 | 52 |
| 4 | Austria | 2 | 26 |
| 5 | Bahrain | 1 | 2 |
| 6 | Pakistan | 1 | 2 |
| 7 | Saudi Arabia | 1 | 2 |
| 8 | Canada | 5 | 136 |
| 9 | Finland | 1 | 8 |
| 10 | New Zealand | 1 | 8 |
| 11 | Sweden | 1 | 8 |
| 12 | France | 2 | 66 |
| 13 | Kazakhstan | 3 | 8 |
| 14 | United states | 4 | 1 |
| 15 | Malaysia | 2 | 34 |
| 16 | Spain | 2 | 11 |
| 17 | Taiwan | 1 | 185 |
| 18 | Iran | 1 | 18 |

In the collaboration network visualization in Figure 9, four clusters have been found among the 13 connected countries. The clusters are separated into four colors (green, yellow, blue, and red). The more citations associated with a country, the bigger that country's node circle. From this, it can be concluded that the UK has the most influential and cited research on composite manufacturing supply chains. Besides the UK, countries like Canada, France, and Germany also have solid international co-authorship with others like Bahrain, Greece and the United States.

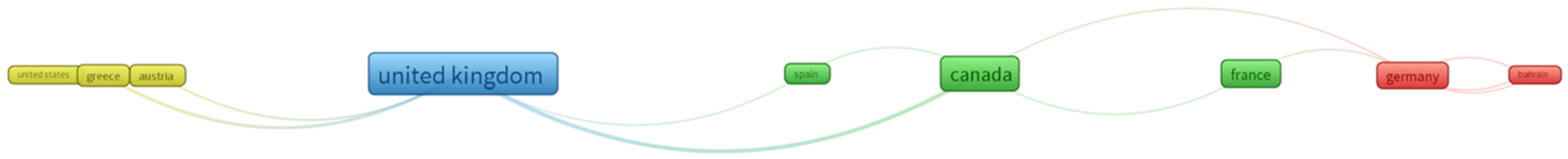

Figure 9: Co-authorship of countries (VosViewer)

## 4. Comprehensive review

This section conducted an in-depth review of the 30 papers selected through the PRISMA framework. Table 4 describes the objectives, types of industry and composite materials, methodology used, and the contribution and outcomes of the articles. Optimization of locations, reducing cost, supply of wastes, and network design are the most dominant goals in composite manufacturing supply chain research (Rentizelas et al. 2022; Shehab et al. 2020; Doustmohammadi and Babazadeh 2020; Vo Dong et al. 2019; Sidelnikov et al. 2021; Shehab et al. 2019; Omair et al. 2022; Trivyza et al. 2022; Rentizelas and Trivyza 2022). Reduction of the supply chain stages or activities is another common objective in this field (Clark and Busch 2002; Mondragon and Mondragon 2017). Moreover, technology and supplier selection are other research purposes among the researchers in their studies (Coronado Mondragon et al. 2016; Hsu et al. 2012; Coronado Mondragon et al. 2017). Furthermore, blockchain, recycling, resiliency, and supply chain behavior have emerged in the objectives analysis among researchers (Mastrocinque et al. 2015; Shuaib et al. 2015; Mondragon et al. 2018; Mondragon et al. 2019; Piri et al. 2018; Meiirbekov et al. 2021; Iyer et al. 2023). In the methodological approach, the most widely used research techniques are case studies and surveys. The Mixed Integer Linear Programming Model (MILP) is commonly applied by different researchers when developing the optimization model. Some researchers also used the system dynamic model, conceptual framework, and empirical study in their strategies, as shown in the table. The most popular composite material supply chain types among the selected articles are Carbon Fiber Reinforced Polymers (CFRP), as shown in Table 4.

Table 4: Summary of selected (n=30) articles.

| **Ref** | **Objectives** | **Types of Composite materials** | **Methodology** | **Major contribution/ outcomes** |
|---|---|---|---|---|
| (Coronado Mondragon et al., 2016) | Investigating variables pertaining to supply chain setup and | Not specific | Survey based | Composites company's respondents rated quality, labor cost reduction, cycle time reduction, and rework |
| | technology choice. | | | highest in technology selection. |
| (Hsu et al., 2012) | Supplier selection of Aluminum composite panel. | Aluminum composite | Empirical study by Decision-making trial and evaluation laboratory (DEMATEL) | Six dimensions are identified (quality of the product (D1), schedule (D2), supply risk (D3), overall cost of product (D4), service (D5), and environmental collaboration (D6)) in respect of supplier selection. |
| (Mastrocinque et al., 2015) | Exploring the composite materials sector to comprehend the supply chain behavior. | Not specific | Conceptual framework | Suggested that the connection between agility, products, greenness, and firms' performance needs to be explored more for the composite materials supply chain. |
| (Shuaib et al., 2015) | End life of composite wastes and benefit of composite materials remanufacturing. | CFRP and Glass fiber reinforced polymer | Sankey diagrams | Demonstrated the potential for cross-sector reuse of recycled composite materials, thereby mitigating the escalating issue of composite waste management. |
| (Mondragon et al., 2015) | Product manufacturing, supplier and technology selection for the configuration of the composite materials supply chain. | Not specific | Case study | Addressed in regard to how supplier selection, lean concepts, and distribution network optimization are becoming more and more important when choosing technologies. |
| (Mahmood et al., 2012) | Investigating green supply chain management in aero composite industries. | Not specific | Survey-based questionnaires and semi-structured interviews | This paper pointed out that optimization of the process and lean manufacturing system are the most common factors for green supply chain practice among the investigated companies. |
| (Clark and Busch, 2002) | Shortening the supply chain for composite parts. | Glass fiber reinforced polymers | Case study | Observed how vertical manufacturing integration helped to cut costs while shortening the supply chain. |
| (Rentizelas et al., 2022) | Supply chain network design problem for blade waste in Europe. | Glass fiber reinforced polymers | Mixed-Integer Linear Programmin g (MILP) model | Pointed out that for sustainable reverse supply network high volume of waste material has a significant role. |
| (Coronado Mondragon et al., 2017) | Factors and challenges of technology selection. | Not specific | Survey based | Out of the eighteen factors related to the technology selection process in composite materials manufacturing supply chains, "supply chain performance" has the greatest number of correlations. |
| (Shehab et al., 2020) | Determining challenges associated with cost modeling. | Carbon Fiber Composites | Secondary data analysis | Identified three issues (technical, supply chain and networking, market challenges) for the recycling of the carbon fiber composites. |
| (Mondragon et al., 2018) | Investigating the applicability of block chain for the supply chain of composite .materials | Carbon Fiber Composites | Secondary data analysis | Exploring the concept of utilizing blockchain technology to provide tamper-proof information sharing between the various supply chain manufacturing stages. |
| (Mondragon et al., 2019) | Feasibility study of block chain technology in carbon fiber and perishable goods supply chain | Carbon Fiber Composites | Case study | Identified that the supply chain enabled by blockchain is useful for the certification of composite material components. |
| (Vaidya and Hopkins Jr, 2021) | Efforts of IACMI on sustainable composite materials manufacturing in | Not specific | Case study | Described the method of producing a composite mask and the adaption policy for the recycling and circular economy. |
| | response to COVID-19 | | | |
| (Piri et al., 2018) | Cause and effect analysis of the three aspects (interconnectivity, adaptability, transformability) of bio-composite production system resiliency. | Bio-composite | System dynamic model | Suggested that interconnectivity, adaptability, and transformability can lead to developing a resilient system. |
| (Mativenga, 2019) | Determining the composite wastes recycling locations | Carbon, glass, resin, and dry fiber | Survey-based and case study. | When the supply chain complexity is reduced the greenhouse gas emission for transportation is also decreased. |
| (Doustmohammadi and Babazadeh, 2020) | Minimizing the total cost of direct and reverse flows in the supply chain | Wood Plastic Composite | MILP model | Closed-loop supply chain model for the composite material is sensitive to the changes in demand and shipping cost. |
| (Mondragon and Mondragon, 2017) | Reducing the number of operations in manufacturing supply chain | Carbon fiber | Case study | Reducing the number of operations from 28 to 10 in the supply chain manufacturing. |
| (Coronado Mondragon et al., 2018) | Investigating the substitution of traditional materials with composite materials to reconfiguration of supply chain, and reduction of activities. | CFRP | Case study | Developed a design process for the adoption of composite materials and end up with a reduction in the number of operations of the supplier network. |
| (Vo Dong et al., 2019) | Optimizing a CFRP waste supply chain by including the cost minimization/ NPV | CFRP | MILP model and case study | Microwave recycling process have several advantages over others based on economic profit and Global warming potential including end-of-life |
| | maximization and Global Warming Potential (GWP). | | | waste acceptability and better recovery of fiber. |
| (Sidelnikov et al., 2021) | Creating a mathematical model to reduce flow from the composite recycling unit and maximize return on investment. | Polymer composites | Mathematical model for the optimal result (nonlinear model) | Determined the optimal volume of composite materials. |
| (Shehab et al., 2019) | Developing a cost modelling framework for recycling CFC materials. | Carbon Fiber Composites | Survey based | Establishing the major cost drivers of mechanical and thermal recycling processes (machines, process time and complexity, materials, labor, and production volumes). |
| (Meiirbekov et al., 2021) | Investigated the uncertainties associated with the recycling of CFRP. | CFRP | Narrative literature review | Found four main categories of uncertainty in the supply chain, including sourcing, demand for recycled goods, and transportation factors. |
| (de los Mozos and López, 2020) | Controlling the short-term logistical operations. | Fiber-reinforced plastics | Case study | Standardized the logistical tasks carried out by the operator, resulting in a decrease in execution time and error. |
| (Badea et al., 2018) | Examine the performance and criterion of wind blades using composite materials while taking supply chain constraints into account. | Carbon fiber, and Fiberglass | Analytical Hierarchy Process (AHP) and Philosophy Drum-Buffer-Rope | Investigated the demand buffer while maintaining the constant of orders at each level of composites material supply chain. |
| (Omair et al., 2022) | Showed the optimal production and | CFRP | A nonlinear mathematical model | Optimal production and outsourcing quantity |
| | outsourcing quantity that will minimize the total cost of the supply chain. | | | minimized the total cost of the supply chain. |
| (Trivyza et al., 2022) | Investigating the cross sectorial circular economy pathways through reverse supply chain network design. | CFRP | MILP model | Optimal network design is more centralized with only two facilities proposed |
| (Rentizelas and Trivyza, 2022) | Reverse supply chain network design for the reusable CFRP car frame. | CFRP | MILP model | Optimal scenario is more centralized with only one facility is required. |
| (Helbig et al., 2016) | Estimating the geo-political supply risk factors. | Carbon fiber | Case study | Identified the geopolitical risk factors. |
| (Iyer et al., 2023) | Decision making in balancing recycled content claims with costs and uncertainties. | Glass fiber | Developing a model for recycled input carryover decisions with supply uncertainty and external sourcing costs. | Recycled content claims are impacted by shorter planning horizon. |
| (Hasan and Wuest, 2022) | Exploring the state of the art of composite materials supply chain. | Not specific | Systematic literature review | Identified the trend of the composite manufacturing supply chain. |

**4.1.Composites supply chain and recycling.**

Most researchers concentrated on the supply chain's recycling industry (n=6) for the composites supply chain, illustrated in Figure 10.

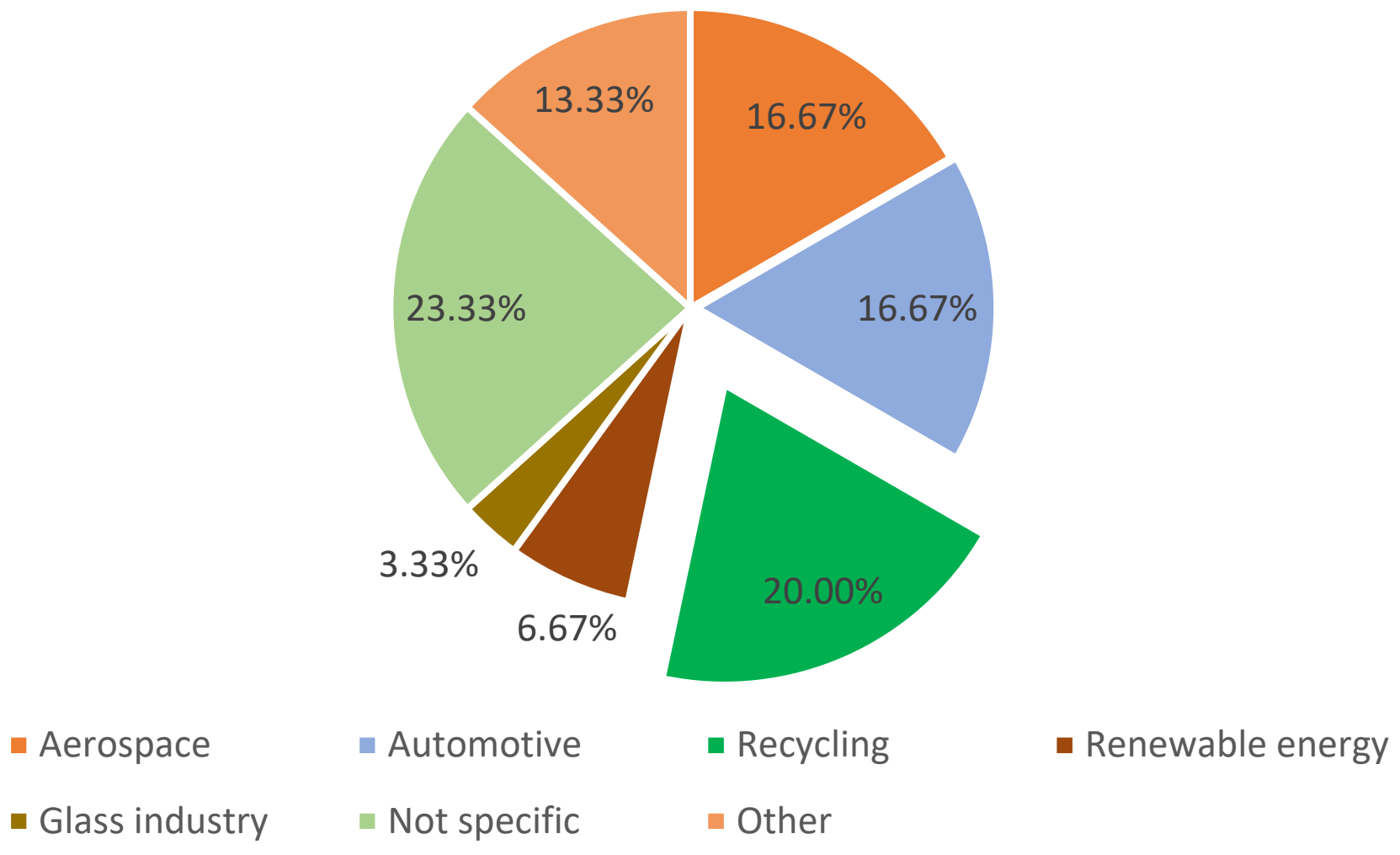

Figure 10: Types of industries in composite materials supply chain

Mastrocinque et al. (2015) urged implementing green practices in the composite manufacturing supply chain by reevaluating the company's supply chain. In 2019, Sultan and Mativenga pointed out that stakeholders have differing views on the waste ownership of composite materials. Partially, the take-back approach has been implemented in Germany, Sweden, and the UK. This is another issue for the circular supply chain of composite materials. Some argue that users who turn the product into waste should be responsible for the product's end life, while others believe the government should take it entirely. In their investigation, one limitation is that they used the K-means algorithm to find optimum clusters for recycling centers. However, the K-means clustering method uses Euclidean distance rather than Mahalanobis distance, which is weighted (Mativenga, 2019).

The trade-off between more extensive facilities with better economies of scale versus decentralized smaller facilities with low transportation costs could be an investigation scenario for the future. Composite waste landfilling is a major environmental problem. For example, only Germany has commercial facilities in Europe for glass fiber-reinforced polymers recycling. The investigation by Rentizelas et al. (2022) is one of the first studies on the supply chain network design problem for end-of-life wind blades. Their analysis predicted that in 2050, the demand for wind blade composites could not be satisfied with only UK facilities. Instead, it will need to partly come from the large facility in Germany. This study is supported by Shehab et al. (2020), who highlighted establishing a market for recycled carbon fibers as an essential point. However, recycled carbon fiber has limitations due to degradation in mechanical and surface properties. Recycling CFRP has a tremendous beneficial environmental impact. Shuaib et al. (2015) described that the energy needed by the pyrolysis method for recycling and remanufacturing CFRP is less than 10% of the total energy associated with manufacturing virgin CFRP materials.

Interestingly, although recycling process efforts are increasing, Vo Dong et al. (2019) surprisingly did not find any public information about a decreasing trend of waste in the Aerospace

CFRP industry. So, this substantial increase scenario lowers the price of recovered carbon fiber due to high production wastes. Determination of recycling sites has been constituted by the closeness or amount of waste. Legislation-related uncertainties like having a potential privilege granted by the government to recycle rather than landfill is also essential, argued by Meiirbekov et al. (2021). In this case, European initiatives have successfully increased the recovery and recycling of post-consumer waste. For instance, Switzerland boasts a 95% recovery rate for post-consumption glass. Manufacturers must take advantage of Green Public Procurement by providing increased recycled content at a given price point, thus cultivating strong demand for recycled content (Iyer et al., 2023). Mahmood et al. (2012) indicated that quality, cost, delivery, and continuous improvement are the desired performance measures for green supply chain management among surveyed companies. They also suggested that the lean manufacturing system is a crucial element of green supply chain practice not only to eliminate waste but also to develop the supplier network, production flow, and less inventory.

**4.2.Supply chain and its drivers**

Six key elements influence the supply chain: Facility, Information, Sourcing, Demand, Pricing, and Transportation (Chopra et al., n.d.). The Figure's 11 shows the number of articles in percentage selected for this study as well as the supply chain drivers associated with them.

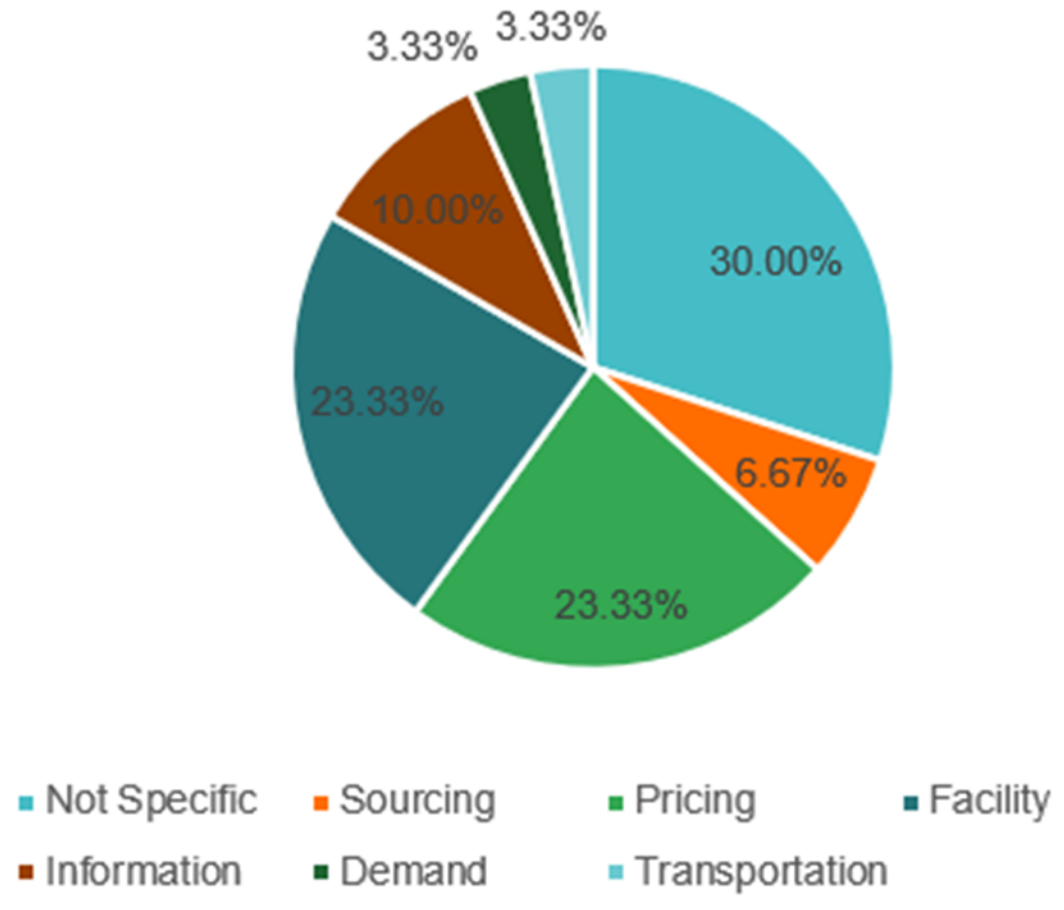

Figure 11: Supply chain drivers

Figure 11 reveals that there has been a gradual rise in pricing related articles in the composites supply chain, although most papers do not fall into a specific supply chain driver category (not specific (n=9)).

An effective supply chain network is essential to meet unanticipated demand while minimizing costs. A study by Trivyza et al. (2022) described that optimal network design is more centralized, with only two facilities proposed in Europe given the low volume of carbon fiber available in 2023. Their model will become more decentralized, with up to five facilities proposed by 2050. Interestingly, their data shows the breakeven price of recycled carbon fiber is competitive against virgin carbon fiber. However, the quality could improve during recycling, and it is more often used in automobiles than aerospace. Substantial evidence found that transported material density helps reduce overall cost but does not have a critical impact. Different cost drivers are used in developing mathematical models. Omair et al. (2022) considered carbon emission which substantially impacts total cost. Surprisingly, their model shows production and inspection costs have less income than total cost. A different study by Sidelnikov et al. (2021) considered the cost of disposal of items that cannot be reprocessed in their model to find the optimal volume of composite materials.

Doustmohammadi and Babazadeh (2020) determined the optimum material flow for the closed-loop supply chain network to minimize total cost. Somewhat surprisingly, changes in transportation and production cost parameters do not greatly affect the total cost. Supplier or vendor selection is crucial for supply chain management. Hsu et al. (2012) depicted that the best vendor of recycled material has good scores on product quality and environmental collaboration. Moreover, they showed supply risk would affect delivery service and schedule, product cost, and environmental collaboration.

Coronado Mondragon et al. (2018) suggested reducing lead time, number of suppliers, and supply chain complexity can positively affect composite material use in manufacturing. For high-tech industries (where flexible and innovative approaches are expected) to maintain stringent standards among manufacturers, authentic information and traceability are essential. Mondragon et al. (2018) highlighted blockchain technology would enhance information flow between composite material supply chain stages with public and private keys. In their follow-up 2019 study, they used a cryptographic hash function for the blockchain application. Besides the information supply chain driver, their survey showed on-time deliveries, supply chain performance, and reduced cycle time are critical for composite raw material producers. They introduced the UK as having industrial base composite expertise. Reducing cycle time significantly influences aerospace, automotive, marine, and motorsport industries (Coronado Mondragon et al., 2017).

### 4.3.KPI’s model for composite materials supply chain.

According to Parmenter (2019), "Key Performance Indicators (KPIs) are used to improve performance dramatically. KPIs concentrate on the vital elements for current and future success. The KPIs were measured subjectively in the comprehensive literature review of selected articles, as shown in the table. It was discovered the KPIs (n=105) used varied widely, as described in Appendix A. Some interesting identified KPIs consist of demand for recycled products, data sharing and confidentiality, quality, cost, delivery, continuous improvement (QCDC), and rate and efficiency of recycling. Social, economic, and environmental are the three pillars of sustainability, a concept gradually emerging from critiques reconciling economic growth to solve ecological and social problems addressed by the United Nations (Purvis et al., 2019). The three circles diagram of sustainability was first introduced by Barbier, as shown in Figure 12.

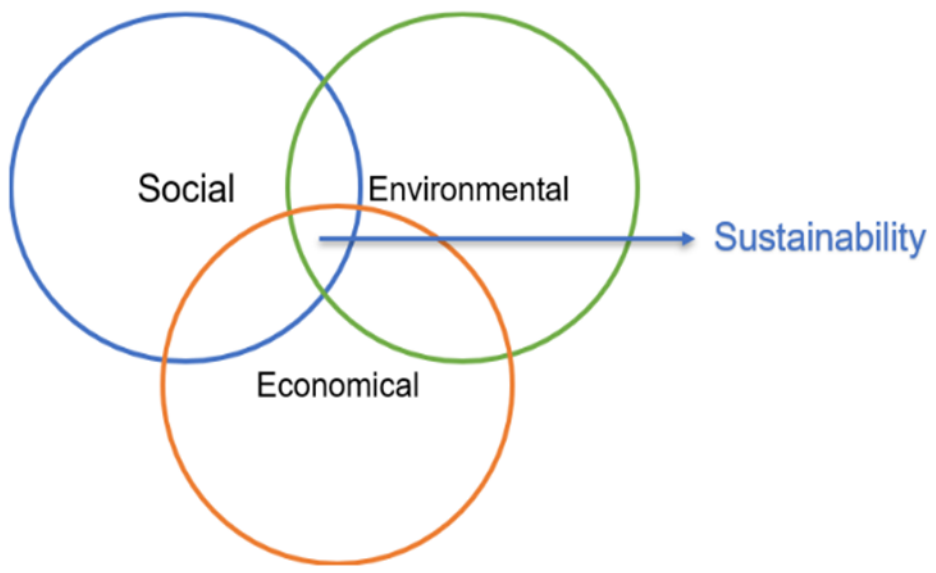

Figure 12: Triple bottom line (TBL) of sustainability (Barrier, 2002).

After merging similar KPIs, they are all modeled with TBL and associated strategic goals. The goal-based KPIs are shown in Table 5, representing the standardized KPIs set found among the selected papers (n=30) on the composite materials supply chain under the TBL umbrella. This KPI model ensures strategic goals advance in a socially responsible, environmentally sustainable, and economically viable way.

Table 5: Standardized KPIs of Composite materials supply chain

| | **KPIs** | | |
|---|---|---|---|
| **Environmental** | • Carbon footprint<br>• Greenhouse gas emissions | • Recycling rate<br>• Waste generation and reduction<br>• Recycled carbon fiber quality | • Energy consumption<br>• Resource-efficient manufacturing<br>• Green practices<br>• Processing facility locations and capacities (when related to environmental impact) |
| Goals | ✓ Reduce carbon footprint | ✓ Enhance waste management and recycling | ✓ Promote resource efficiency and sustainable manufacturing |
| **Social** | • Customer satisfaction<br>• Demand fulfillment<br>• Order fill rate<br>• Customer service level | • Supplier selection<br>• Delivery schedule (impacts work-life balance and fair working hours) | • Quality of the product (affects customer safety and satisfaction)<br>• Data sharing and confidentiality<br>• Supply chain transparency and traceability<br>• Supply chain collaboration (related to ethical practices and corporate social responsibility) |
| Goals | ✓ Improve customer satisfaction and engagement | ✓ Ensure fair labor practices and supplier ethics | ✓ Foster a culture of safety, quality, and transparency |
| **Economical** | • Profit/loss | • Technology selection | • Overall cost of product |
| | • Cost reduction | • Supply chain performance | • Annual revenues from avoided waste disposal |
| | • Return flow volume optimization | • Complexity and links | • Market value and volume |
| | • Material recovery rate | • Agile and lean practices | • Cost savings |
| | • Processing time | • Economic efficiency | • NPV of recycling fiber |
| | • Recycling efficiency | • Supplier reliability | • Global warming potentials (economic impacts) |
| | • Supply uncertainties | • Efficiency in procurement | • Minimization of total cost |
| | • Demand of recycled products | • Reliability of supply chain network | • Transportation costs |
| | • Inventory cost | • Manufacturing agility | • Legislation (economic impacts) |
| | • Inventory management | • Supply chain resiliency and efficiency | • Raw material availability |
| | • Cost efficiency | • Facility utilization | • Market responsiveness |
| | • Recycling cost | • Total cost of the reverse supply chain | |
| | • Inventory holding cost | | |
| | • Inspection cost | | |
| | • Production rate | | |
| | • New material launch rate | | |
| Goals | ✓ Maximize economic efficiency and profitability | ✓ Optimize supply chain performance and reliability | ✓ Ensure sustainable cost management and investment returns |

### 4.4.Challenges and future works

Several questions still need answering to improve and maintain sustainable composites supply chain. For example, shortening the supply chain operation is tricky. Clark and Busch (2002) highlighted matching scales in vertically integrated processes is critical. To ensure efficient vertical integration, integrated operations' scales must be comparable to realize significant economic benefits. Determining recycled carbon fiber products' value can be complex due to various factors like origin, composition, manufacturing time, and costs. Mondragon et al. (2019)

described that composite material production varies by industry, complicating supply and demand balancing. Composites contain carbon, glass, and aramid fibers and resins like epoxy and polyester, each with distinct, complex properties posing supply challenges. The demand for environmentally friendly materials and processes intensifies the challenge of finding sustainable raw materials, manufacturing methods, and recycling solutions for composite waste (Mondragon et al., 2019). When deciding on adopting composite materials, there are multiple factors to consider, including cost, feasibility, and operational benefits, which can complicate the composite materials supply chain (Mondragon and Mondragon, 2017). This has been further supported by Vaidya and Hopkins Jr (2021), who described manufacturers must consider recycled materials' costs and availability when making claims to drive demand. Obtaining usable recycled materials from nearby sources can be complex and unpredictable.

Another area for improvement in the composite materials supply chain is data sharing. Shehab et al. (2020) said industry members' reluctance to share information about their processes, pricing, and competitive details can significantly impact cost model accuracy. Recognizing transparency and collaboration's value in achieving reliable cost estimates is essential. Working together and sharing information can ensure cost models are as accurate as possible, ultimately benefiting the industry. Mativenga (2019) also argued making well-informed waste management strategy decisions is challenging without precise data. Mondragon et al. (2018) described the composite materials industry's lack of standardized materials and procedures makes implementing a uniform blockchain system in the supply chain difficult. Compliance with regulations and standards is non-negotiable for industries like medical and aerospace. So, data security and privacy are crucial in decentralized blockchains, especially for sensitive manufacturing information. Again, the constantly changing legislation and policies surrounding waste disposal and recycling make effective implementation difficult for companies. Companies need to adapt to these changes while continuing to comply with regulations (Mativenga, 2019). Different certification agencies have different planning horizons for verifying claims, which can affect recycled content claims' viability due to interaction with external sourcing costs (Iyer et al., 2023).

The overall challenges described in the selected articles (n=30) are shown in Appendix B. Similar challenges are grouped together (e.g., "Data availability and reliability" and "Real-time data integration" combined as "Data availability and reliability"). The challenges are then categorized under three fundamental decision-making levels - strategic, tactical, and operational - a strategy adopted from the 2018 review paper by Barbosa-Póvoa et al. (2018). The cluster analysis of challenges is illustrated in Figure 13. Strategic planning is long-term, tactical planning is shorter-term concerned with demand and inventory planning, and operational planning relates to current tasks needing execution (Miller 2012; Ivanov 2010).

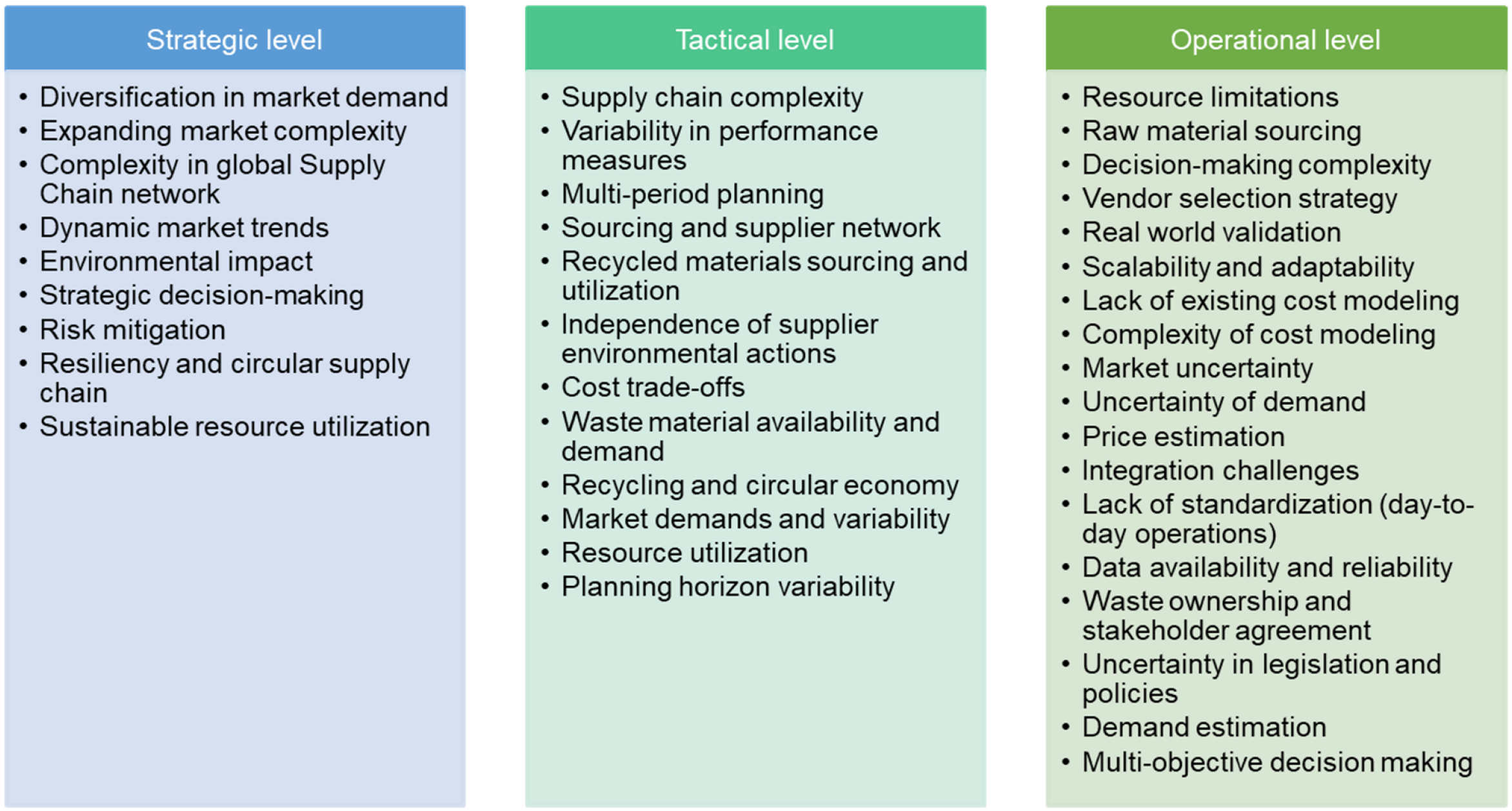

Figure 13: Cluster analysis of challenges based on decision levels.

This categorization allows the composite materials supply chain to prioritize challenges based on their impact on long-term goals, medium-term objectives, and day-to-day operations. This chart also helps decision-makers efficiently allocate resources to each challenge based on the appropriate time frame. Moreover, this categorization of significant challenges provides communication among the different supply chain echelons for a common understanding of each challenge type and its nature. Furthermore, it fosters continuous improvement efforts among different tiers of the supply chain by reviewing performance against different time-scaled objectives.

This study further highlighted some future research directions:

- Establishing the composites recycled market to improve material flow in the composite materials supply chain.
- Identifying a common data-sharing platform through blockchain technology integration.
- Effective way of collaboration and selection of the suppliers of composite materials supply chain.
- Establishing legislation for recycled content claims of composite materials.
- Analyzing sustainable resource utilization and the circular supply chain.
- Domino or ripple effect analysis to investigate supply chain resiliency.
- Improving supply chain performance, including the drivers of supply chain.
- Evaluating the composite materials supply chain based on market demand.
- Enhancing material and information flow to better understand complex multi-echelon, multi-objective composite material supply chains.
- Optimizing the composite materials supply chain transportation and distribution network.

- Synchronizing forward and reverse supply chains to avoid operational disruptions.

## 5. Conclusion

This study investigates the overview of the sustainable and green supply chain for composite materials. This study makes directions on how to close the existing research gaps and connect sustainability and green initiatives in composite materials supply chains in the future. There has been a growing number of articles in this field recently. Although the bibliometric and comprehensive review analyzed the field while maintaining the PRISMA systematic literature review framework, this study has limitations in what could be accomplished. There are areas for improvement in how results are presented and structured. Expanding keywords in the search terminology may enable a more exhaustive field review. At the same time, a larger paper pool could require more advanced bibliometric tools to overcome large, clean, combined dataset difficulties. Further content analysis of specific manuscripts may elucidate future directions and current research progress.

The bibliometric analysis shows an ongoing positive trend among researchers in this field, where most research has diverse themes needing more collaboration. It has also been observed that a small number of researchers have generated a significant number of influential results. Some prominent trending research fields are recycling, optimization, cost modeling, composite waste reverse supply chain, and MILP-based studies, where most research articles are either survey-based or case studies. The comprehensive analysis of the 30 selected papers indicated various challenges, and focusing on the identified KPIs is necessary to improve the current and future composite materials supply chain. As this research field is still maturing, there are ample opportunities for scholars to advance it and its subfields in the coming days.

## Declarations

### Funding

No funding was received for this research.

### Conflict of Interest/ Competing interests

The authors declare no conflict of interest or competing interests.

### Availability of Data/Materials

Not available for publication.

### Code availability

Not available for publication.

### Ethics approval

The authors declare the article was constructed respecting all ethical conditions of publication.

**Consent to participate.**

All authors participated in the preparation of the article. Therefore, the authors allow their names in the article.

**Consent for publication**

The authors allow for publication.

**CRediT authorship contribution statement**

Md Rabiul Hasan: Conceptualization, Investigation, Methodology, Data curation, Formal analysis, Writing – original draft, Visualization, Writing – review & editing. Muztoba Ahmed Khan: Visualization, Writing – review and editing, Validation, Supervision. Thorsten Wuest: Writing – review and editing, Validation, Supervision, Project administration.

**Declaration of Competing Interest**

The authors declare that they have no known competing financial interests or personal relationships that could have appeared to influence the work reported in this paper.

## Appendix A: KPI's of composite materials supply chain

| Ref. | KPI's | Ref. | KPI's |
|---|---|---|---|
| (Coronado Mondragon et al., 2016) | 1. Technology selection<br>2. Supply chain performance | (Doustmohammadi and Babazadeh, 2020) | 1. Facility utilization<br>2. Waste reduction<br>3. Customer satisfaction<br>4. Recycling rate |
| | | | 5. Minimization of total cost<br>6. Demand fulfillment |
| (Mondragon et al., 2015) | 1. Supplier selection<br>2. Technology selection<br>3. Complexity and links | (Mondragon and Mondragon, 2017) | 1. Lead time<br>2. Cost<br>3. Number of operations |
| (Mastrocinque et al., 2015) | 1. Green practices<br>2. Agile and lean practices | (Coronado Mondragon et al., 2018) | 1. Supplier locations<br>2. Number of operations |
| (Shuaib et al., 2015) | 1. Resource efficient manufacturing<br>2.Environmental footprints | (Vo Dong et al., 2019) | 1. NPV of recycling fiber<br>2. Global warming potentials |
| (Hsu et al., 2012) | 1. Quality of the product<br>2. Delivery schedule<br>3. Supply risk<br>4. Overall cost of product | (Sidelnikov et al., 2021) | 1. cost reduction<br>2. Return flow volume optimization<br>3. Recyclability rate |
| (Mahmood et al., 2012) | 1. Effectiveness of Lean Manufacturing System (LMS) and tools<br>2. Quality, Cost, Delivery, and Continuous Improvement (QCDC) | (Shehab et al., 2019) | 1. Material recovery rate<br>2. Processing time<br>3. Recycling efficiency |
| (Clark and Busch, 2002) | 1. Supply Chain Length<br>2. Economic Efficiency | (Meiirbekov et al., 2021) | 1. Supply uncertainties<br>2. Demand of recycled products<br>3. Transportation<br>4. Legislation<br>5. Inconsistent market price of recycled CF products |
| (Rentizelas et al., 2022) | 1. Carbon emissions per processing stage/ Carbon footprint<br>2. Processing facility locations and capacities | (de los Mozos and López, 2020) | 1. Order fill rate<br>2. Customer service level<br>3. Inventory cost<br>4. Order cycle time and delivery reliability |
| | 3. Annual revenues from avoided waste disposal<br>4. Profit/loss | | |
| (Shehab et al., 2020) | 1. Supplier reliability<br>2. Energy consumption<br>3. Recycled carbon fiber quality<br>4. Data Sharing and Confidentiality | (Badea et al., 2018) | 1. Raw material availability<br>2. Lead time<br>3. Inventory management<br>4. Cost efficiency<br>5. Quality control<br>6. Supply chain collaboration |
| (Mondragon et al., 2018) | 1. Reduction of lead times<br>2. Efficiency in procurement<br>3. Supply chain transparency and traceability | (Omair et al., 2022) | 1. Carbon emission<br>2. Recycling cost<br>3. Waste reduction<br>4. Inventory holding cost<br>5. Inspection cost<br>6. Production rate |
| (Mondragon et al., 2019) | 1. Real-time visibility<br>2. Certification process efficiency | (Trivyza et al., 2022) | 1. Facility utilization<br>2. Carbon emission<br>3. Total cost of the reverse supply chain<br>4. Facility utilization |
| (Vaidya and Hopkins Jr, 2021) | 1. Market value and volume<br>2. Sustainability metrics<br>3. Resilience and adaptability<br>4. Composite material sales<br>5. Waste generation and reduction | (Rentizelas and Trivyza, 2022) | 1. Total system cost<br>2. Carbon emission cost<br>3. Capacity utilization of the facilities |
| (Piri et al., 2018) | 1. Reliability of supply chain network<br>2. Manufacturing agility<br>3. Compliance with regulations<br>4.Supply chain resilience<br>5. Market responsiveness | (Helbig et al., 2016) | 1. Production concentration<br>2. Disclosed information about the production and import volumes<br>3. Geopolitical supply risk factor |
| | 6 New materials launch rate | | |
| (Mativenga, 2019) | 1. Recycling rate<br>2. Processing capacity utilization<br>3.Carbon footprint and greenhouse gas Emissions<br>4.Cost savings | (Iyer et al., 2023) | 1. Recycled content claim Impact<br>2. External sourcing costs<br>3. Planning horizon<br>4. Demand benefit<br>5.Raw material cost |
| (Coronado Mondragon et al., 2017) | 1.Supply chain performance and technology selection rating | (Hasan and Wuest, 2022) | 1. Collaboration among researchers<br>2. Innovation in recycling and end of life processes<br>3. Supply chain resiliency and efficiency<br>4. Market demand for recycled materials |

Appendix B: Major challenges of composite materials supply chain

| Reference | Challenges | Reference | Challenges |
|---|---|---|---|
| (Mondragon et al., 2016) | 1. Diversification in market demand<br>2.Complexity in supply chain<br>3. Variability in performance measures<br>4. Resource limitations | (Mativenga, 2019) | 1. Supply chain complexity<br>2. Data availability and reliability<br>3. Multi-objective decision making<br>4. Waste ownership and stakeholder agreement<br>5. Uncertainty in legislation and policies |
| (Mondragon et al., 2015) | 1. Lack of standardization<br>2. Expanding market complexity | (Doustmohammadi and Babazadeh, 2020) | 1. Raw material sourcing<br>2. Closed loop collection and recycling<br>3. Distribution and collection center placement<br>4. Demand estimation<br>5. Environmental impact<br>6 Multi-period planning |
| (Mastrocinque et al., 2015) | 1. Complexity in global Supply Chain network<br>2. Dynamic market trends | (Mondragon and Mondragon, 2017) | 1. Lack of standardization<br>2. Supply chain complexity<br>3. Sourcing and supplier network<br>4 Decision-making complexity |
| (Shuaib et al., 2015) | 1. Resource efficiency<br>2. Decision-making challenges | (Coronado Mondragon et al., 2018) | 1. Integration challenges<br>2. Standardization and processes |
| (Hsu et al., 2012) | 1. Vendor selection strategy<br>2. Comprehensive and integrated approach<br>3. Interdependencies and feedback<br>4. Recycled materials sourcing and utilization | (Vo Dong et al., 2019) | 1. Price estimation |
| (Mahmood et al., 2012) | 1. Supply chain dynamics and environmental pressure<br>2. Independence of supplier environmental actions | (Sidelnikov et al., 2021) | 1. Real world validation<br>2. Scalability and adaptability |
| (Clark and Busch, 2002) | 1. Sourcing complexity<br>2. Cost trade-offs<br>3. Matching scale of operations | (Shehab et al., 2019) | 1. Lack of existing cost modelling<br>2. Complexity of cost modeling |
| (Rentizelas et al., 2022) | 1.Market uncertainty<br>2. Waste material availability and demand<br>3. Reverse supply network optimization<br>4. Feasibility of circular economy pathways | (Meiirbekov et al., 2021) | 1. Market acceptance of the recycled carbon fiber products<br>2. Quality of recovered composites<br>4. Uncertainty of demand |
| (Shehab et al., 2020) | 1. Waste supplier co-operation<br>2. Lack of industry disclosure<br>3. Market establishment<br>4. Variation in raw material quality<br>5. Product pricing and labeling | (de los Mozos and López, 2020) | 1. High demand uncertainty<br>2. Real time data integration<br>3. inventory management<br>4. Operational disruptions |
| (Mondragon et al., 2018) | 1. Lack of standardization<br>2. Energy-intensive manufacturing<br>3. Complex supply chain<br>4. Regulatory Compliance<br>5 Data security and privacy<br>6.Change management and adoption | (Badea et al., 2018) | 1. Restriction of supply chain collaboration<br>2. Need for buffer stock<br>3. Negative impact of raw material availability |
| (Mondragon et al., 2019) | 1. Integration and data sharing<br>2.Data accuracy<br>3. Sustainability<br>4. Market demands and variability<br>5 Complexity of raw materials | (Omair et al., 2022) | 1. Stochastic demand<br>2. Waste management<br>3. Multi-echelon supply chain management<br>4. Policy of circular economy |
| (Vaidya and Hopkins Jr, 2021) | 1. Sourcing of usable recycled input<br>2. Balancing demand-side benefits<br>3. Recycling and circular economy<br>4. Embodied energy and environmental Impact | (Trivyza et al., 2022) | 1. Multi-objective and time frame optimization problem<br>2. Carbon emission and sustainability |
| (Piri et al., 2018) | 1. Diverse market needs<br>2.Resource utilization<br>3.Strategic decision-making<br>4.Risk mitigation<br>5. Variability in raw materials | (Rentizelas and Trivyza, 2022) | 1. Holistic approach of the circular economy<br>2. Stochastic optimization |
| (Helbig et al., 2016) | 1. Lack of domestic production<br>2. Fluctuations in prices and availability<br>3. Dependencies on imports | (Iyer et al., 2023) | 1. Multi-manufacturer competition<br>2. Diverse standards and product features<br>3. Planning horizon variability |
| (Hasan and Wuest, 2022) | 1.Resiliency and circular supply chain<br>2. Improve the efficiency<br>3. Sourcing strategies<br>4. Optimization of the cost and distribution network<br>5. Cascading failure or domino-effect analysis<br>6. Sustainable resource utilization | (Coronado Mondragon et al., 2017) | 1.Limitation of resources<br>2. Variability of performance measures. |